\documentclass[12pt]{article}
\usepackage{theorem, amsfonts,amsmath,graphicx}
\usepackage{setspace}        
\nonstopmode
\def\CA{{\cal A}}

\def\CF{{\cal F}}

\def\pa{\partial}

\newcommand\COMP{\hbox{C\kern -.58em {\raise .54ex
\hbox{$\scriptscriptstyle |$}}\kern-.55em {\raise .53ex
\hbox{$\scriptscriptstyle |$}} }}

\newcommand\MM{\hbox{I\kern-.2em\hbox{M}}}
\newcommand\NN{\hbox{I\kern-.2em\hbox{N}}}
\newcommand\RR{\hbox{I\kern-.2em\hbox{R}}}
\newcommand\sRR{{\it \hbox{I\kern-.2em\hbox{R}}}}
\newcommand\QQ{\hbox{I\kern-.53em\hbox{Q}}}
\newcommand\PP{\hbox{I\kern -.2em\hbox{P}}}
\newcommand\EE{\hbox{I\kern-.2em\hbox{E}}}
\newcommand\ZZ{{{\rm Z}\kern-.28em{\rm Z}}}
\newcommand\II{{{\rm I}\kern-.28em{\rm I}}}

\newcommand\qed{\hfill$\sqcap\kern-8.0pt\hbox{$\sqcup$}$}

\newcommand\be{\begin{equation}}
\newcommand\ee{\end{equation}}

\title{Indifference pricing and hedging in stochastic volatility models}

\author{M. R. Grasselli and 
 T. R. Hurd  \\
Dept. of Mathematics and Statistics\\
McMaster University\\ Hamilton, ON, L8S 4K1}
\date{April 23, 2004}
\begin{document}
\maketitle

\begin{abstract} 
We apply the concepts of utility based pricing and hedging of derivatives in stochastic volatility markets and introduce a new
class of ``reciprocal affine" models for which the indifference price and optimal hedge portfolio for pure volatility claims are efficiently computable.  We obtain a general formula for the market price of volatility risk in these models and calculate it explicitly for the case of an exponential utility. 
\end{abstract}

\noindent
{\bf Key words:} volatility risk, exponential utility, reciprocal affine models, Davis price, incomplete markets.

\newpage
\section{Introduction}

The purpose of this paper is to provide a concrete application of {\em utility based
pricing and hedging} for markets with untraded risk factors.  For incomplete markets of this type, they provide an economically justified framework for pricing derivatives where the traditional paradigms of replication and risk neutral valuation fail to produce any identifiable price. For instance, this framework is well suited for pricing in insurance markets (see Young and Zariphopoulou 2002) \nocite{YouZar02}
as well as market models with stochastic volatility. 

The method is based on solutions of the so called {\em optimal hedging problem} (see Cvitanic, Schachermayer and Wang 2001; Owen 2002), \nocite{CviSchWan01, Owen02}
which is a generalization of Merton's optimal investment problem, whose solution has a distinguished history
including key works such as Karatzas, Lehoczky, Shreve  and Xu (1991)
\nocite{KaLeShXu91}, \nocite{KraSch99} and \nocite{Scha01a}
Kramkov and  Schachermayer (1999) and Schachermayer (2001). Characteristically, these papers present abstract existence and uniqueness results for the optimal portfolios, but do not address the question of actually constructing or systematically approximating them given specific market models. However, because of the variety and richness offered by general incomplete markets, a full understanding of the appropriateness of this pricing concept can only achieved
by applying it to an increasing number of concrete examples.

The valuation of derivatives in markets with stochastic volatility typically depends on the choice of
an appropriate dynamics for the volatility process $\sigma_t$, the specification of a functional form for the volatility risk premium and the use of asymptotic filtering techniques to estimate the model parameters and the current level of the unobservable process $\sigma_t$ in terms of the time series of observed stock prices $S_t$ 
(see Fouque, Papanicolaou and Sircar  2000; Gr\"{u}nbichler  and Longstaff 1996; Heston 1993; Stein and Stein 1991). 
\nocite{FoPaSi00,GruLon96,Hest93,SteSte91}
While we do not address estimation and calibration issues, we show how the utility based framework offers insight for the choice of stochastic volatility models and provides a specification of volatility risk premia.

We consider contingent claims on  two factor stochastic volatility models of the form
\begin{eqnarray}
d\bar S_t&=& \bar S_t[\mu(t,Y_t)dt + \sigma(t,Y_t) dW^1_t], \nonumber\\
dY_t&=&a(t,Y_t)dt + b(t,Y_t)[\rho dW^1_t +\sqrt{1-\rho^2}dW^2_t],
\label{stocvolSDE}
\end{eqnarray}
with initial values $\bar S_0,Y_0\geq 0$, deterministic functions $\mu, \sigma, a, b$ satisfying the regularity and growth conditions necessary to ensure the existence and uniqueness of a solution to \eqref{stocvolSDE}, independent one dimensional $P$--Brownian motions $W^1_t$ and $W^2_t$, and finally constant correlation $|\rho| < 1$. 

In section 2, we review the concepts of {\em certainty equivalent} and {\em indifference price} of a claim as introduced by  Hodges and Neuberger (1989) \nocite{HodNeu89}
 (according to Becherer 2001) \nocite{Bech01}) 
for a general utility function and then specialize to an exponential utility $U(x)=-e^{-\gamma x}$ with risk aversion parameter $\gamma > 0$. We then show how the utility based approach uniquely specifies the volatility risk premium associated with each contingent claim $B(S_T,Y_T)$, as well as a claim independent market price for volatility risk which is related to the solution of Merton's problem for a stochastic volatility market.

We use the theory of dynamic programming in section 3 to derive the nonlinear PDE satisfied 
by the certainty equivalent of a general contingent claim on the market specified by \eqref{stocvolSDE}.
When the contingent claim is a pure volatility claim (such as in Brenner and Galai 1989; Whaley 1993),
\nocite{BreGal89,Whal93}
we will see in section 4 how a change of variables leads to a linear second order PDE, for which a Feynman-Kac representation of the solution is available. 

In section 5, we explore the striking resemblance of this formula to the framework of interest rate derivatives. We introduce  
the {\em reciprocal affine} stochastic volatility model, where the ``spot rate" $R_t \propto \sigma^{-2}_t$ is taken to be a CIR process. 

In section 6 we use a transform analysis technique developed by Duffie, Pan and Singleton (2000)
\nocite{DuPaSi00} 
to obtain tractable expressions for the indifference price and 
the holdings in both the hedging and Merton portfolios for stochastic volatility. The choice of a reciprocal affine model for the squared volatility process fully reveals its intrinsic elegance in this section when we compute its market price for volatility risk explicitly. In section 7 we benefit from the efficiently implementable formulas derived previously and present numerical results for volatility put options and the difference of volatility call options, discussing several qualitative features of the indifference prices and the hedging strategies for such instruments. In particular, for $\gamma$ not too large, we observe clearly that pricing and hedging are approximately independent of the risk aversion.

\section{Utility based pricing} 

\subsection{General setup}

The central definition is that of an {\em optimal hedging portfolio}. It corresponds to the strategy followed by an investor with initial wealth $x$ who, when faced with a (discounted) financial liability $B$ maturing at a future time $T$, tries to solve the stochastic control problem
\be
u(x)=\sup_{\CA} E\left[U\left(X_T - B\right)|X_0=x\right],
\label{optimalhedge}
\ee
where $X_T$ is the discounted terminal wealth obtained when investing $H_t\bar S_t$ dollars on the risky asset and $\eta_tC_t$ dollars in a riskless cash account with value $C_t$ initialized at $C_0=1$ and governed by
\be
\label{BDE}
dC_t=r_t C_t dt.
\ee
The utility function $U:\RR\to\RR\cup\{-\infty\}$ is assumed to be a concave, strictly increasing and differentiable function and the domain of optimization  ${\cal A}$ is restricted to self--financing portfolios, 
that is, to wealth processes satisfying
\begin{equation}
C_tX_t:=H_t\bar S_t+\eta_tC_t=x+ \int_0^t H_ud\bar S_u+\int_0^t\eta_udC_u.
\label{intwealth}
\end{equation}
In addition to the self--financing condition, economic reasoning imposes further restrictions on the class 
${\cal A}$ of admissible portfolios. For instance, one might require that the wealth process $X_t$ be uniformly bounded from below, in order to rule out ``doubling strategies" (following Harrison and Pliska 1981).
\nocite{HarPli81}
Finally, the liability 
$B$ is assumed to be a random variable of the form $B=B(S_T,Y_T)$, for some Borel function 
$B:\RR^2_+\to \RR$.
 
It follows from It\^{o}'s formula and the definitions \eqref{stocvolSDE} and \eqref{BDE} that the market model 
in terms of the discounted price process $S_t=\bar S_t/C_t$ is 
\begin{eqnarray}
dS_t&=& S_t[(\mu(t,Y_t)-r)dt + \sigma(t,Y_t) dW^1_t],  \nonumber\\
dY_t&=&a(t,Y_t)dt + b(t,Y_t)[\rho dW^1_t +\sqrt{1-\rho^2}dW^2_t].
\label{stocvolSDE2}
\end{eqnarray}
We have taken $r_t=r$ to be constant for simplicity; the case of either a deterministic rate $r_t=r(t)$ or a stochastic 
rate of the form $r_t=r(t,S_t,Y_t)$ are immediate but artificial generalizations. The really relevant case of stochastic rates driven by a different stochastic factor, possibly correlated with $S_t$ and $Y_t$, is not pertinent to the scope of this paper.

Using the self-financing condition \eqref{intwealth}, we immediately obtain that the discounted 
wealth process satisfies 
\begin{equation}
\label{xSDE}
dX_t = H_tdS_t=H_tS_t[(\mu(t,Y_t)-r)dt+\sigma(t,Y_t)dW^1_t],
\end{equation}
so that the only relevant control in \eqref{optimalhedge} is $H_t$, with the holdings in the cash account being determined by  
$\eta_t=X_t-H_tS_t$.

For Markovian markets such as \eqref{stocvolSDE2} and claims of the form $B=B(S_T,Y_T)$,
we can embed the optimal hedging problem \eqref{optimalhedge} into the larger 
class of optimization problems defined by 
\begin{equation}
\label{dynamic}
u(t,x,s,y)=\sup_{H\in{\cal A}_t}E_{t,s,y}[U(X_T-B(S_T,Y_T))|X_t=x], 
\end{equation}
for $t\in(0,T)$,
where $x\in\RR$ denotes some arbitrary level of wealth, ${\cal A}_t$ denotes admissible
portfolios starting at time $t$ and $E_{t,s,y}[\cdot]$ denotes expectation with respect to the joint probability law at time $t$ of the processes $S_u,Y_u$ satisfying \eqref{stocvolSDE2} for 
$u\geq t$, with initial condition $S_t=s$ and $Y_t=y$.

\subsection{The certainty equivalent and the indifference price}

Suppose that \eqref{dynamic} has an optimizer $H^B_t$, that is, assume that 
\[u(t,x,s,y)=E_{t,s,y}[U(x+(H^B\cdot S)_t^T-B(S_T,Y_T))],\]
where we use the notation 
\[(H\cdot S)_t^T:=\int_t^T H_u dS_u.\]
Define the {\em certainty equivalent} for the claim $B$ at time $t$ as the process $c^B_t=c^B(t,x,s,y)$ satisfying the equation
\begin{equation}
U(x-c^B_t))=E_{t,s,y}[U(x+(H^B\cdot S)_t^T-B(S_T,Y_T))].
\label{cert}
\end{equation}
In other words, it corresponds to the amount which when subtracted from the wealth $x$ at time $t$ produces the same deterministic value for the utility as the optimal expected utility of terminal wealth starting at $x$ and facing 
the claim $B$ at terminal time $T$. If we set $B=0$, then the optimal hedging problem becomes the Merton optimal investment problem and we denote the certainty equivalent by $c^0_t=c^0(t,x,s,y)$. 

An agent with utility $U$ and wealth $x$ at time $t\in(0,T)$ will charge a premium for issuing a liability $B$ 
maturing at $T$ . The {\em indifference price} for the claim $B$ is defined to be the premium that makes the agent indifferent between  making the deal or not, that is, the unique solution $\pi^B=\pi^B(t,x,s,y)$ (if it exists) to the equation
\begin{equation}
\sup_{H\in{\cal A}_t}E_{t,s,y}[U(x+(H\cdot S)_t^T]=\sup_{H\in{\cal A}_t}E_{t,s,y}[U(x+\pi^B+(H\cdot S)_t^T-B(S_T,Y_T)].
\end{equation}
From the definition of the certainty equivalent, we see that this equation is equivalent to
\begin{equation}
U(x-c^0_t)=U(x+\pi^B-c^B_t),
\end{equation}
so that the indifference price is given by
\begin{equation}
\pi^B(t,x,s,y)=c^B(t,x+\pi^B(t,x,s,y),s,y)-c^0(t,x,s,y).
\end{equation}
The {\em Davis price} (see Davis 1997) \nocite{Davi97} 
of the contingent claim $B$ is defined as the linearization of the indifference price, that is,
\begin{equation}
\pi_{Davis}^B:=\left.\frac{d(\pi^{\varepsilon B})}{d\varepsilon}\right|_{\varepsilon=0}.
\end{equation}

As we see, establishing the existence of an indifference price is tantamount to solving both the optimal hedging problem for the claim $B$ and the Merton investment problem. As such, it depends on the choice of the utility function and the appropriate notion of
admissible portfolios, as well as integrability conditions on the claim $B$ itself. 

For the rest of this paper we concentrate 
on exponential utilities of the form
\begin{equation}
\label{exp}
U(x)=-e^{-\gamma x},
\end{equation}
where $\gamma>0$ is the risk--aversion parameter. The advantage of using an exponential utility is that  we can factorize the value function $u(t,x,s,y)$ in \eqref{dynamic} as
\begin{eqnarray}
u(t,x,s,y)&=&\sup_{H\in{\cal A}_t}E_{t,s,y}\left[-e^{-\gamma\left(x+(H\cdot S)_t^T-B(S_T,Y_T)\right)}\right] \nonumber \\
&=& -e^{-\gamma x} \inf_{H\in{\cal A}_t}E_{t,s,y}\left[e^{-\gamma\left((H\cdot S)_t^T-B(S_T,Y_t)\right)}\right] \nonumber \\
&=:&  U(x)v(t,s,y).
\label{factor}
\end{eqnarray}

To define the domain of optimization, let ${\cal M}^a(S)$ and ${\cal M}^e(S)$ denote respectively the sets of absolutely continuous and equivalent (local) martingale measures for $S$ and let ${\cal M}^f(S)$ denote the set of measures $Q \in {\cal M}^a(S)$ with finite relative entropy with respect to $P$. For absence of arbitrage, we assume that ${\cal M}^e\cap {\cal M}^f\neq \emptyset $. Now let $L(S)$ denote the set of predictable $S$--integrable processes. For definiteness, we take the set of admissible portfolios to be
\begin{equation}
{\cal A}=\{H \in L(S) : (H \cdot S)_t \mbox{ is a  Q--martingale for all Q } \in {\cal M}_f\},
\end{equation}
but note that alternative choices for ${\cal A}$ also yield the existence and uniqueness of optimal hedging portfolio for the exponential utility. It follows from the abstract convex duality solution of the optimal hedging problem that, for claims $B$ satisfying certain integrability conditions, the certainty equivalent 
derived using the exponential utility exists as a well defined process (see Becherer 2001, 2004; Delbaen, Grandits, Rheinl{\"a}nder, Samperi, Schweizer and Stricker 2002; Kabanov and Stricker 2002; Owen 2002).
\nocite{DGRSSS02, KabStr02, Owen02, Bech01,Bech04}

It follows directly from \eqref{cert} that the certainty equivalent is wealth independent and given by
\begin{equation}
\label{logv}
c^B(t,s,y)=\frac{1}{\gamma} \log v(t,s,y).
\end{equation}
Analogously, setting $B=0$ gives the certainty equivalent for the Merton problem with exponential utility as
\begin{equation}
c^0(t,s,y)=\frac{1}{\gamma} \log v^0(t,s,y),
\end{equation}
where $v^0(t,s,y)= \displaystyle{\inf_{H\in{\cal A}_t}}E_{t,s,y}\left[e^{-\gamma(H\cdot S)_t^T}\right]$. Thus the indifference price process for the claim $B$ obtained from an exponential utility is given by
\begin{equation}
\pi^B(t,s,y)=c^B(t,s,y)-c^0(t,s,y)=\frac{1}{\gamma} \log\frac{v(t,s,y)}{v^0(t,s,y)}.
\label{piexp}
\end{equation}

\subsection{Martingale Pricing}

The solution of the hedging problem \eqref{optimalhedge} through convex duality involves finding a martingale 
measure $Q^B$ minimizing the expected value of the Legendre transform of the utility function $U$ over the set ${\cal M}^a(S)$ of absolutely continuous martingale measures for $S$. The optimal wealth $X^B_T=x+(H^B\cdot S)_T$ and the optimal martingale measure $Q^B$ are then related by the fundamental equation 
\begin{equation}
\label{fund}
U^\prime(X^B_T-B)=\xi \frac{dQ^B}{dP},
\end{equation}
where $\xi=u^\prime(x)$ (see Owen 2002). Letting 
\begin{equation}
\Lambda^B_t:=E_t\left[\frac{dQ^B}{dP}\right]=\frac{1}{\xi}E_t[U^\prime(X^B_T-B)]
\end{equation}
be the density process for the measure $Q^B$, we define the {\em utility based price of risk} associated with the claim $B$ as the vector process $\lambda^B_t=((\lambda^B)^1_t,(\lambda^B)^2_t)$ satisfying
\begin{equation}
\frac{d\Lambda^B_t}{\Lambda^B_t}=-[(\lambda^B)^1_t dW^1_t + (\lambda^B)^2_t dW^2_t].
\end{equation}

For the case of exponential utility this reduces to 
\begin{eqnarray}
\Lambda^B_t&=&-\frac{\gamma}{\xi}E_t[U(X^B_t+(H^B\cdot S)_t^T-B] \nonumber \\
&=& -\frac{\gamma}{\xi} U(X^B_t-c^B_t) \nonumber \\
&=& \frac{e^{-\gamma(X^B_t-c_t^B)}}{e^{-\gamma(x-c^B_0)}},
\end{eqnarray}
where we have used the definition of the certainty equivalent process $c^B_t$ and the fact that $\xi=u^\prime(x)=U^\prime(x-c^B_0)$. Applying It\^{o}'s lemma to the martingale $\Lambda^B_t$ above, we easily obtain that 
\begin{eqnarray}
(\lambda^B)^1_t&=&\gamma[(H^B-\partial_sc^B)\sigma S-b\rho \partial_y c^B ],  \label{lambda1}\\
(\lambda^B)^2_t&=&-\gamma b \partial_y  c^B \sqrt{1-\rho^2}  \label{lambda2}.
\end{eqnarray}

Both the optimal martingale measure $Q^B$ and the utility based price of risk $\lambda^B_t$ are specifically related to the 
claim $B$ and therefore do not constitute a direct generalization of the paradigm of pricing by expectation with respect to a risk  adjusted measure valid for all claims. For instance, as we will later observe in our numerical examples, this is reflected in the fact that the indifference price $\pi^B$ is not linear in the claim $B$, with the obvious effect that $Q^B$ fails to be a pricing measure even for multiples of $B$, let alone for other unrelated claims. For the Davis price, however, matters are slightly more familiar, as the following argument shows.

By differentiating the identity
\[U(x-c^{\varepsilon B}_t)=E_t[U(x+(H^{\varepsilon B}\cdot S)_t^T-\varepsilon B)]\]
at $\varepsilon=0$ we obtain Davis's formula 
\begin{eqnarray}
\pi^B_{Davis}&=&\left.\frac{dc^{\varepsilon B}_t}{d\varepsilon}\right|_{\varepsilon=0} \nonumber \\
&=& \frac{E_t[U^\prime(x+(H^0\cdot S)_t^T)B]}{U^\prime(x-c^0_t)}\nonumber \\
&=& E_t^{Q}[B],
\end{eqnarray}
where we have used \eqref{fund} for Merton's problem. We see that the Davis price can be calculated
as an expectation with respect to the risk adjusted measure $Q$ obtained as the dual solution for the optimal investment problem according to the utility function $U$. This measure induces a {\em utility based market price of risk} valid for all claims, which in the case of an exponential utility is obtained from \eqref{lambda1},\eqref{lambda2} by setting $B=0$.

\section{The PDE  for the certainty equivalent}

Consider the continuous Markovian market specified by \eqref{stocvolSDE2}, an exponential utility function of the form \eqref{exp} and a claim of the form $B=B(S_T,Y_T)$ satisfying appropriate integrability conditions. By the dynamic programming principle, the value function $u(t,x,s,y)$ satisfies the Hamilton--Jacobi--Bellman equation
\begin{eqnarray}
\pa_t u &+&\frac{1}{2}s^2\sigma^2 \partial^2_{ss}u+ b\rho s\sigma \partial^2_{sy}u+\frac{1}{2}b^2 \partial^2_{yy} u
 +s(\mu-r) \partial_s u+a \partial_y u \\
&+&\max_h\left\{\frac12 h^2s^2\sigma^2 \partial^2_{xx} u+b\rho h s \sigma \partial^2_{xy} u + hs^2 \sigma^2 \partial^2_{xs}  u +hs(\mu-r)\partial_x u\right\} =0, \nonumber
\end{eqnarray}
with boundary condition $u(T,z)=-e^{-\gamma (x-B(s,y))}$ and the optimal portfolio process is given by 
\begin{equation}
H^B_t=h^B(t,X_t,S_t,Y_t),
\end{equation}
where $h^B$ is the optimizer of the expression above. Direct substitution of \eqref{factor} into this HJB problem leads to an optimizer of the form $H^B_t=h^B(t,s,y)$, with 
\be
h^B(t,s,y)=\frac{1}{\gamma}\frac{\partial_s v}{v}+\frac{b\rho}{\gamma s \sigma}\frac{\partial_y v}{v}+\frac{(\mu-r)}{\gamma s \sigma^2}.  
\label{Hopt}
\ee
The partial differential equation satisfied by the optimal function $v(t,s,y)$ is then  
\begin{eqnarray}
\partial_t v&+&\frac12\left(s^2 \sigma^2 \partial^2_{ss} v + 2 b\rho s\sigma \partial^2_{ys}  v + b^2 \partial^2_{yy} v\right)
+\left[a-\frac{b\rho(\mu-r)}{\sigma}\right] \partial_y v   \nonumber \\
&&-\frac{1}{2}\left[\frac{1}{v}(b\rho \partial_y v+s \sigma \partial_s v)^2+\frac{(\mu-r)^2}{\sigma^2} v\right] = 0,
\label{vPDE}
\end{eqnarray}
subject to the terminal condition $v(T,s,y)=e^{\gamma B(s,y)}$.

From (\ref{logv}), we find that the certainty equivalent process $c^B(t,s,y)$ is a solution to the 
partial differential equation
\begin{eqnarray}
\partial_t c^B&+&  \frac{1}{2}(s^2\sigma^2\partial^2_{ss} c^B+2s\sigma b\rho \, \partial^2_{sy} c^B + b^2 \partial^2_{yy} c^B)+\left[a-\frac{ b \rho (\mu-r)}{\sigma}\right]\partial_y c^B
\nonumber\\ 
&&-\frac{(\mu-r)^2}{2\gamma \sigma^2}+\frac{\gamma}{2}b^2(1-\rho^2)(\partial_y  c^B)^2=0,
\label{cPDE}
\end{eqnarray}
with terminal condition $c^B(T,s,y)=B(s,y)$. The nonlinearity $(\partial_y v)^2/v$ for the value function has been transformed  into the quadratic nonlinearity $(\partial_y c^B)^2$ for the certainty equivalent, while the nonlinearity associated with 
$(\partial_s v)^2/v$  is
absent. The partial differential equation satisfied by $c^0(t,s,y)$, the certainty equivalent for Merton's problem, is identical to
\eqref{cPDE},  but with the terminal condition $c^0(T,s,y)=0$. 

Using (\ref{Hopt}), the optimal portfolio can be obtained in terms of the certainty equivalent process by
\be
h^B(t,s,y)=\partial_sc ^B+\frac{b\rho}{s\sigma}\partial_y c^B+\frac{(\mu -r)}{\gamma s \sigma^2}.  
\label{opthedge}
\ee
From \eqref{lambda1},\eqref{lambda2} we obtain that the exponential utility based price of risk associated with the claim $B$ is
\begin{eqnarray}
\label{market}
(\lambda^B)^1_t &=& \frac{(\mu(t,Y_t)-r)}{\sigma(t,Y_t)}, \\
(\lambda^B)^2_t &=& -\gamma b \partial_y c^B \sqrt{1-\rho^2},
\end{eqnarray}
with the desirable property that the dependence on the claim $B$ occurs only in its second component, while the component 
related to the first Brownian motion is formally identical to the market price of risk of a complete market. Setting $B=0$ in the equations above gives the market price of volatility risk induced by the exponential utility through the solution of Merton's optimal investment problem. 

\section{Volatility claims}

The equations from the previous section simplify considerably when the claim is independent of the process $S_t$. Because of the exponential nature of 
the SDE for $S_t$ in \eqref{stocvolSDE2}, we obtain that for pure volatility claims of the form $B=B(Y_T)$, the equation for
the certainty equivalent $c^B_t=c^B(t,y)$ is reduced to
\begin{equation}
\partial_t c^B+  \left[a-\frac{b \rho (\mu-r)}{\sigma}\right]\partial_y c^B 
+\frac{1}{2}b^2 \partial^2_{yy} c^B -\frac{(\mu-r)^2}{2\gamma \sigma^2}+\frac{\gamma}{2}b^2(1-\rho^2)(\partial_y c^B)^2=0,
\label{cvol}
\end{equation}
subject to the terminal condition $c^B(T,y)=B(y)$, whereas the optimal hedging portfolio is given by 
\be
h^B(t,s,y)=\frac{1}{s}\left[\frac{b\rho}{\sigma}\partial_y c^B+\frac{(\mu-r)}{\gamma  \sigma^2}\right].  
\label{opthedge1}
\ee

Following  Zariphopoulou (2001)
\nocite{Zari01} we now use the transformation
\begin{equation}
\label{cb}
c^B(t,y)=\frac{1}{\gamma(1-\rho^2)}\log f(t,y),
\end{equation}
to reduce \eqref{cvol}  to the linear parabolic final value problem
\begin{eqnarray}
\partial_t f+  \frac{1}{2}b^2 \partial^2_{yy} f+\left[a-\frac{b \rho (\mu-r)}{\sigma}\right]\partial_y f 
 -\frac{(1-\rho^2)(\mu-r)^2}{2 \sigma^2}f  &=& 0, \\
f(T,y)&=&e^{\gamma(1-\rho^2) B(y)} . \nonumber
\end{eqnarray}

Under the appropriate growth and boundedness assumptions on the coefficient functions $\mu,\sigma,a$ and $b$, we can use the Feynman--Kac formula to represent the solution to the problem above as
\begin{equation}
f(t,y)=\widetilde E _{t,y}\left[e^{-\int_t^T R(s,Y_s)ds} e^{\gamma(1-\rho^2)B(Y_T)}\right],
\label{FK}
\end{equation}
where we define
\begin{equation}
R(t,y)=\frac{(1-\rho^2)(\mu(t,y)-r)^2}{2\sigma(t,y)^2},
\label{Req}
\end{equation}
and $\widetilde E_{t,y}[\cdot]$ denotes the expectation with respect to the probability law at time $s=t$ of the solution to 
\begin{eqnarray}
dY_s&=&\left[a-\frac{b (\mu-r)\rho}{\sigma}\right]ds  +b\left[\rho d\widetilde W^1_s+\sqrt{1-\rho^2}d\widetilde W^2_s\right], \nonumber \\
Y_t &=& y
\label{modY}
\end{eqnarray}
for a pair of independent one dimensional $\widetilde P$--Brownian motions $\widetilde W^1_t,\widetilde W^2_t$, for a probability measure $\widetilde P$ on $(\Omega,\CF,(\CF_t)_{t\in[0,T]})$. If we further required $S$ to be a $\widetilde P$ martingale, the comparison with \eqref{stocvolSDE2} leads to the identification
\begin{eqnarray}
d\widetilde W^1_t &=& dW^1_t+\widetilde\lambda^1_t dt \nonumber \\
d\widetilde W^2_t &=& dW^2_t,
\label{tilde}
\end{eqnarray}
where 
\be
\widetilde \lambda^1_t=\frac{\mu(t,Y_t)-r}{\sigma(t,Y_t)}.
\ee

The martingale measure $\widetilde P$ is a computational device in order to 
obtain the function $f(t,y)$ through \eqref{FK}, and should not be confused with either the optimal martingale measure $Q^B$ associated with the claim $B$ or the risk adjusted measure $Q$ giving rise to the Davis price as explained in section 2. 

\section{Reciprocal affine models}

The first factor in the integrand of the Feynman-Kac representation \eqref{FK} is reminiscent of 
a stochastic discount factor with a
``spot rate" $R_t$ implicitly related to our original non-traded factor $Y_t$ through \eqref{Req}. Following this analogy, we see that when $B=0$, \eqref{FK} is the formal equivalent of the price of a zero coupon bond. If in addition equation \eqref{Req} can be inverted, then the second factor
in \eqref{FK} becomes a claim on the spot rate $R_t$ and we can use the technology developed for interest rate derivatives
in order to price it.

To explore this analogy with interest rates further, for the remainder of the paper we take $\mu$ and $r$ to be constants and  $\sigma(t,Y_t)=\sqrt{Y_t}$, so that
\eqref{Req} becomes
\begin{equation}
R_t=R(t,Y_t)=\frac{(1-\rho^2)(\mu-r)^2}{2Y_t}.
\end{equation}

Affine models form a well-studied class of interest rate models (see Duffie, Pan and Singleton 2000),
\nocite{DuPaSi00}
often leading to analytic expressions for quantities such as bond prices. We can carry the results from these models to our problem by hypothesizing that $R_t$ follows an affine process and thus $Y_t$ is a reciprocal affine process. We illustrate the idea in the specific case of the Cox, Ingersoll and Ross (1985) 
\nocite{CoxIngRos85} 
model. Since our calculations are going to take place under the measure $\widetilde P$, we specify the dynamics for $R_t$ as 
\begin{equation}
\label{CIR}
dR_t=\widetilde\alpha(\widetilde\kappa-R_t)dt + \beta\sqrt{R_t} \left[\rho d\widetilde W^1_t +\sqrt{1-\rho^2}d\widetilde W^2_t \right],
\end{equation}
for constants $\widetilde\alpha,\widetilde\kappa,\beta>0$ with $4\widetilde\alpha\widetilde\kappa>\beta^2$. It follows from 
\eqref{tilde} that the dynamics of $R_t$ under the economic measure $P$ is
\begin{equation}
\label{CIR2}
dR_t=\alpha(\kappa-R_t)dt + \beta\sqrt{R_t} \left[\rho dW^1_t +\sqrt{1-\rho^2}dW^2_t \right],
\end{equation}
where $\alpha=\left(\widetilde\alpha-\beta\rho\sqrt{\frac{2}{1-\rho^2}}\right)$ and $\alpha\kappa=\widetilde\alpha\widetilde\kappa$. We then obtain from the It\^{o} formula that
\begin{eqnarray}
\label{YCIR}
dY_t &=& -\frac{(1-\rho^2)(\mu-r)^2}{2 R_t^2}dR_t + \frac{(1-\rho^2)(\mu-r)^2}{2R_t^3}dR_t^2 \\
&=& \left[\frac{2(\beta^2-\alpha\kappa )}{(1-\rho^2)(\mu-r)^2}Y^2_t+\alpha Y_t\right]dt \nonumber \\
&& -\left(\frac{2}{1-\rho^2}\right)^{1/2}\frac{\beta Y_t^{3/2}}{(\mu-r)}
\left[\rho dW^1_t +\sqrt{1-\rho^2}dW^2_t \right] \nonumber,
\end{eqnarray}
which by comparison with \eqref{stocvolSDE2} leads to  
\begin{eqnarray}
a(t,Y_t)&=& \alpha Y_t + \frac{2(\beta^2-\alpha\kappa)}{(1-\rho^2)(\mu-r)^2}Y^2_t, \\
b(t,Y_t)&=& -\left(\frac{2}{1-\rho^2}\right)^{1/2}\frac{\beta}{(\mu-r)} Y^{3/2}.
\end{eqnarray}

The reciprocal of a CIR process, being positive and mean reverting, is a reasonable model for stochastic volatility. Yet the expressions above show that our formulation differs from popular choices for stochastic volatility, such as the Heston (1993) model,
\nocite{Hest93}
for which $Y_t$ itself follows a CIR process, or the Stein and Stein (1991) model,
\nocite{SteSte91}
where $\sigma_t$ follows an arithmetic 
Ornstein--Uhlenbeck process.

\section{Pricing and hedging formulas}
The indifference price and Davis price of a volatility claim under the reciprocal affine models of the previous section both require computation of expressions of the form
\begin{equation}
I:=\widetilde E_t\left[e^{-\int^T_t R_s ds}g(R_T)\right],
\end{equation}
for functions $g:\mathbb{R}^+\to\mathbb{R}$. Provided its Fourier transform is well defined and invertible, we can 
express $g$ as 
\begin{equation}
g(R)=\frac1{2\pi}\int^{\infty}_{-\infty} e^{-iuR}\hat g(u) du,
\end{equation}
where
\begin{equation}
\hat g(u)=\int_{-\infty}^{\infty} e^{iuR}g(R)dR.
\end{equation}
Using Fubini's theorem to interchange the order of integration, we have
\begin{equation}
I=I(R_t,t,T)= \frac1{2\pi}\int^{\infty}_{-\infty} \Psi (u) \hat g(u) du,
\end{equation}
where $\Psi$ can be computed following Duffie, Pan and Singleton (2000) as
\nocite{DuPaSi00}
\begin{eqnarray}
\Psi(u)=\Psi(u,R_t,t,T)&:=& \widetilde E_t\left[e^{-\int^T_t R_s ds}e^{-iuR_T} \right] \nonumber\\
&=&\exp[M(u,t,T)+N(u,t,T)R_t].
\end{eqnarray}
Here 
\begin{eqnarray}
N(u)=N(u,t,T)&=&\frac{(b_2+iu)b_1-(b_1+iu)b_2 e^{\Delta(t-T)}}{(b_2+iu)-(b_1+iu) e^{\Delta(t-T)}}, \nonumber\\
M(u)=M(u,t,T)&=&\frac{-2\alpha\kappa}{\beta^2}\log{\left(\frac{b_2+iu}{b_2-N}\right)}+\alpha\kappa b_1(t-T),
\end{eqnarray} 
with $b_2>b_1$ being the two roots of $x^2-\frac{2\widetilde\alpha}{\beta^2}x-\frac2{\beta^2}$ and $\Delta=\sqrt{\widetilde\alpha^2+2\beta^2}$.

Setting $g(R_T)=e^{\gamma(1-\rho^2)B(R_T)}$, we obtain from \eqref{piexp}, \eqref{cb} and \eqref{FK} that  the indifference price of the volatility claim $B=B(R_T)$ is simply
\begin{eqnarray}
\pi^B&=&\frac{\delta}{\gamma}\log\left[\frac{\widetilde E_t\left[e^{-\int_t^T R_sds} e^{\gamma(1-\rho^2)B(R_T)}\right]}
{E_t^Q[e^{-\int_t^T R_sds}]}\right]  \nonumber \\
&=& \frac{1}{\gamma(1-\rho^2)} \log\left[\frac{I(R_t,t,T)}{\Psi(0,R_t,t,T)}\right].
\label{pib}
\end{eqnarray}

From \eqref{opthedge1}, \eqref{cb} and \eqref{FK}, the number of shares of stock to be held in order to optimally 
hedge against the claim $B$ is
\begin{eqnarray}
h^B(t,y)&=&\frac{1}{\gamma s}\left[\frac{b\rho}{\gamma(1-\rho^2)\sqrt{y}}\frac{\partial \log I}{\partial y}+\frac{(\mu-r)}{\gamma y}\right] \nonumber \\
&=& \frac{1}{\gamma s}\frac{(\mu-r)}{y}\left[\frac{\beta\rho}{\sqrt{2(1-\rho^2)}}
\frac{\int_{-\infty}^{\infty}\Psi(u)N(u)\hat g(u)du}{\int_{-\infty}^{\infty}\Psi(u)\hat g(u)du}+1\right],
\label{hb}
\end{eqnarray}
whereas the number of shares held in the Merton portfolio is
\begin{eqnarray}
h^0(t,y)&=&\frac{1}{s}\left[\frac{b\rho}{\gamma(1-\rho^2)\sqrt{y}}\frac{\partial \log \Psi(0)}{\partial y}+\frac{(\mu-r)}{\gamma y}\right] \nonumber \\
&=& \frac{1}{\gamma s}\frac{(\mu-r)}{y}\left[\frac{\beta\rho}{\sqrt{2(1-\rho^2)}}N(0)+1\right].
\label{h0}
\end{eqnarray}

Finally, using \eqref{market} with $B=0$, the reader is invited to perform an amusing calculation leading to the conclusion that the market price of volatility risk induced by an exponential utility in the reciprocal affine model is
\begin{equation}
\left(\begin{array}{c}
\lambda^1_t \\ \lambda^2_t
\end{array}\right)=
\left(\begin{array}{c}
1 \\ \frac{\beta}{\Delta\sqrt{2}}(1-e^{\Delta(t-T)})
\end{array}\right)
\frac{(\mu-r)}{\sqrt{Y_t}}.
\end{equation}
This remarkable functional form for the market price of volatility risk has several pleasant implications. We first observe that in the regime of a sufficiently remote time to maturity both of its components are constant multiples of the usual market price of risk obtained for complete markets with a time 
varying volatility $\sigma_t=\sqrt{Y_t}$. Secondly, for typical parameter values, the first component is about two orders of magnitude bigger than the second one, implying that most of the volatility risk is already accounted by the dependence of the model on the first Brownian motion $W^1_t$. Finally, it is now easy to verify that our process $Y_t$ is the reciprocal of a CIR process under {\em both} the economic measure $P$ and the exponential risk adjusted measure $Q$, a consistency result which is far from obviously true for other stochastic volatility models.

\section{Numerical results}

Numerical computation of prices and hedge amounts is now highly efficient using the fast Fourier transform. We found that with a discrete Fourier lattice of $2^{12}$ points, we could compute hundreds of prices and hedging amounts per second with accuracy better than 1\% on a desktop PC. This speed enabled us to easily survey dependences of prices and hedging in all parameter values.

We illustrate the range of possibilities for model parameters fixed  at reasonable values: $\rho=0.5, \alpha=5, \beta=0.04, \kappa=0.001, \mu=0.04, r=0.02$ and initial squared volatility ranging in the interval $[0, 0.5]$. With these parameters the squared volatility process has a mean reversion time of approximately two months and an equilibrium distribution with expected value approximately 40\%.

Figures 1, 2 and 3 present plots of the indifference price of  a volatility put option with payoff $(K-\sigma^2_T)^+$ with strike $K=0.15$ for risk aversions 
$\log_2\gamma$ in the range $[-5,5]$ and maturities $T$ in the range $[0.1,1]$ years. Figure 4 shows the dollar amount of the hedge in excess of the Merton strategy for the volatility put as a function of $Y_0$ and $T$ for $\gamma=1$. These experiments showed in particular that there is no significant dependence of these values on $\gamma$ in the range $[0,10]$. We can also observe how, for longer maturities, the mean reverting character of $Y_t$ leads prices which are practically independent of $Y_0$.

The present model has the somewhat surprising property that unbounded claims such as forward contracts and call options have an expected utility of negative infinity. For that reason, we restrict ourselves to claims for which the payoff is bounded from above, such as the difference of calls. Figure 5 shows  the indifference price for a difference of volatility calls $(\sigma^2_T-K_1)^+-(\sigma^2_T-K_2)^+$ with $K_1=0.15, K_2=0.3$ for  maturities $T$ in the range $[0.1,1]$ years and risk aversion $\gamma=1$.

In figure 6 we select one sample path for the squared volatility process $Y_t=\sigma^2_t$ and the indifference price for the volatility put with $K=0.15, T=1, \gamma=1$. We then show the difference between the amount invested in the stock according to the hedging portfolio and  the amount invested in the stock following Merton's strategy. This plot shows clearly that long before maturity, the hedger can ignore the contingent claim, but that they become increasingly anxious closer to maturity. Further experiments show again that these values are rather insensitive to the choice of $\gamma $ in the range $[0,10]$.

\section{Conclusions}

The main goal of this paper has been to find a financially viable class of stochastic volatility models and related derivative securities which lead to tractable forms for the utility based indifference price and the corresponding optimal hedging strategies. As is well known, the indifference price is wealth independent only for exponential utility, so we have adopted this as the preferred utility. The indifference price and its linearization the Davis price are motivated by fundamental economic principles applicable to general incomplete markets, which makes these questions of  theoretical interest and of use to practitioners using models with either stochastic volatility or jumps.

We have found that ``close to closed form solutions''  can be derived very elegantly for contingent claims which depend solely on  spot volatility at maturity. A class of stochastic volatility models for which this is possible was found which is as rich as but inequivalent to the original Heston model. Contingent claims on spot volatility have been studied before in the literature, but always with an implied assumption about the market price of risk. Our treatment is the first to use the market price of risk derived from first principles via the utility function.   A number of the interesting properties of volatility claims have been explored 
by other authors (such as in Gr\"{u}nbichler  and Longstaff 1996), and these effects can be easily observed in our formulation.

Our formulas can be efficiently computed using the fast Fourier transform: the resulting algorithm yields hundreds of prices per second on a desktop PC. This speed gives us the capability to explore the indifference price, optimal hedging and sensitivities and the market price of risk for a complete range of financially relevant model parameters including risk aversion, mean rates of return, mean reversion rates, correlation parameter, mean volatility and volatility of volatility. We have shown a number of the possibilities in this paper, but many detailed properties remain to be explored in further papers.

Amongst the important questions still to be addressed we mention two. The first is to quantify the performance of our optimally  hedged portfolios relative to unhedged portfolios for a variety of volatility claims. Since the volatility claims we analyze are traditionally considered to be ``unhedgable'', it is of interest to use quantities such as value--at--risk to measure portfolio performance. We expect to observe performance which is directly related to the correlation between volatility and stock price. A second question is to compare the exponential utility based market price of risk we compute with the forms assumed in the literature by other researchers. We hope to show that our computed market price of risk  gives reasonable consistency with observed option price data.

\end{document}